\documentclass{amsart}

\usepackage[english]{babel}
\usepackage[all]{xy}
\usepackage{latexsym}
\usepackage{amssymb}
\usepackage{amsmath}
\usepackage{amsthm}
\usepackage{amscd}
\usepackage{a4wide}

\theoremstyle{plain}
\newtheorem{pro}{\sc Proposition}
\newtheorem{thm}[pro]{\sc Theorem}

\newtheorem*{thmintro}{\sc Th\'eor\`eme}
\theoremstyle{definition}
\newtheorem*{dei}{\sc Definition}

\newtheorem*{Rq}{\sc Remark}
\newtheorem*{exem}{\sc Examples}

\newtheorem*{sketch}{\sc Sketch of the proof}

\newcommand{\B}{\mathcal{B}}

\newcommand{\oPo}{\bar{\Po}}

\newcommand{\BLi}{\mathcal{B}i\mathcal{L}ie}
\newcommand{\IBi}{\mathcal{I}nf\mathcal{B}i}
\newcommand{\G}{\mathcal{G}}
\newcommand{\N}{\mathcal{N}}

\newcommand{\Sy}{\mathbb{S}}

\newcommand{\Po}{\mathcal{P}}

\newcommand{\F}{\mathcal{F}}
\newcommand{\ac}{\scriptstyle \textrm{!`}}

\newcommand{\Q}{\mathcal{Q}}

\newcommand{\fin}{\ \hfill \square}
\newcommand{\Co}{\mathcal{C}}

\begin{document}
Projet de Note aux C. R. Acad. Sci. Paris

Algebra / \emph{Alg\`ebre}

\vskip20pt

\begin{center}
{\huge \ \bf Koszul Duality for PROPs}\\
\vspace{15pt}

{\large Bruno Vallette\footnote{Institut de Recherche
Math\'ematique Avanc\'ee, Universit\'e Louis Pasteur et CNRS,  7,
Rue Ren\'e Descartes, 67084 Strasbourg Cedex, France.\\
Email : \texttt{vallette@math.u-strasbg.fr}}}
\end{center}

\vskip20pt

\hrule \vspace{15pt}

\textbf{Abstract}\\

The notion of PROP models the operations with multiple inputs and multiple outputs,
acting
on some algebraic structures like the bialgebras or the Lie bialgebras.
We prove a Koszul duality theory for PROPs generalizing the one for
associative algebras and for operads.\\

\textbf{R\'esum\'e}\\

\textbf{Dualit\'e de Koszul des PROPs. } La notion de PROP
mod\'elise les op\'erations \`a plusieurs entr\'ees et plusieurs
sorties, agissant sur certaines structures alg\'ebriques comme les
big\`ebres et les big\`ebres de Lie. Nous montrons une th\'eorie
de dualit\'e de Koszul pour les PROPs qui g\'en\'eralise celle des
alg\`ebres associatives et
des op\'erades.\\

\hrule

\vskip25pt

\subsection*{Version fran\c caise abr\'eg\'ee}$\ $\\

Suivant J.-P. Serre dans \cite{Serre}, on regroupe sous le terme
de \emph{g\`ebre} diff\'erentes structures alg\'ebriques comme les
alg\`ebres, les cog\`ebres et les big\`ebres.\\

L'ensemble $\Po(m,\, n)$ des op\'erations \`a $n$ entr\'ees et $m$ sorties agissant
sur un
certain type de g\`ebres est un module \`a gauche sur le groupe sym\'etrique
$\Sy_m$ et \`a droite sur $\Sy_n$.

On appelle \emph{$\Sy$-bimodule} toute collection $(\Po(m,\, n)
)_{m,\, n \in \mathbb{N}^*}$ de tels modules. Nous d\'efinissons
un produit $\boxtimes$ dans la cat\'egorie des $\Sy$-bimodules qui
repr\'esente les compositions d'op\'erations \`a plusieurs
entr\'ees et plusieurs sorties. Ce produit est bas\'e sur les
graphes dirig\'es (\emph{cf.}
Figure~\ref{figrapheniveau}).\\

On d\'efinit un \emph{PROP} comme un $\Sy$-bimodule muni d'une
composition $\Po \boxtimes \Po \xrightarrow{\mu} \Po$. On donne
les exemples du PROP $\BLi$ des big\`ebres de Lie (\emph{cf.}
\cite{D}), du PROP $\BLi_0$ des big\`ebres de Lie combinatoires
(\emph{cf.} \cite{C}) et du PROP $\IBi$ des big\`ebres de Hopf
infinit\'esimales (\emph{cf.} \cite{A}). On appelle
\emph{$\Po$-g\`ebre}, tout module sur le PROP $\Po$. On retrouve
les d\'efinitions des g\`ebres classiques. Par exemple, une
$\BLi$-g\`ebre est exactement une big\`ebre de
Lie.\\

Nous \'etendons les d\'efinitions de bar et cobar constructions
des alg\`ebres et des op\'erades aux PROPs, et nous
g\'en\'eralisons les lemmes de comparaison  de B. Fresse
\cite{Fresse} aux PROPs. Remarquons que les d\'emonstrations
op\'eradiques ne sont pas reconductibles ici, car ces derni\`eres
reposent sur les propri\'et\'es combinatoires des arbres.\\

A  partir d'un PROP gradu\'e par un poids, par exemple quadratique
c'est-\`a-dire d\'efini par des g\'en\'erateurs et des relations
quadratiques, on construit une \emph{coPROP} dual $\Po^{\ac}$ et
un \emph{complexe de Koszul} $(\Po^{\ac}\boxtimes \Po,\, d_K)$.\\

Le principal th\'eor\`eme de cette th\'eorie
est le suivant :

\begin{thmintro}
Soit $\Po$ est un PROP diff\'erentiel augment\'e gradu\'e par un
poids (par exemple, un PROP quadratique), les propositions
suivantes sont \'equivalentes
\begin{enumerate}
\item le complexe de Koszul $(\Po^{\ac}\boxtimes \Po,\, d_K)$ est
acyclique,

\item  la cobar construction sur le coPROP dual $\Po^{\ac}$
fournit une r\'esolution du PROP $\Po$ :
$$ \bar{\B}^c(\Po^{\ac})\longrightarrow \Po.$$
\end{enumerate}
\end{thmintro}

Dans ce cas, la r\'esolution obtenue est le \emph{mod\`ele
minimal} de $\Po$ et elle permet de d\'efinir la notion de
\emph{$\Po$-g\`ebre \`a homotopie pr\`es}.\\

Nous montrons que le PROP $\BLi$ des big\`ebres de Lie (\emph{cf.}
\cite{D}), le PROP $\BLi_0$ des big\`ebres de Lie combinatoires
(\emph{cf.} \cite{C}) et le PROP $\IBi$ des big\`ebres de Hopf
infinit\'esimales (\emph{cf.} \cite{A}) sont des PROPs de Koszul.
Ce qui permet de donner les d\'efinitions de big\`ebres de Lie,
big\`ebres de Lie combinatoires et big\`ebres de Hopf
infini\'esimales \`a homotopie pr\`es.

\vskip10pt \hrule \vskip30pt

\section*{\bf Introduction}

The Koszul duality theory for algebras, proved by S. Priddy in
\cite{P} has been generalized to the operads by V. Ginzburg and M.
M. Kapranov in \cite{GK}.

An operad models the operations acting on a certain type of
algebras (associative, commutative and Lie algebras for instance).
Since these operations have multiple inputs but only one output,
their compositions can be represented by trees. This theory has a
lot of applications. It gives the minimal model of an operad
$\Po$, the notion of $\Po$-algebras up to homotopy and a natural
homology theory for the $\Po$-algebras.

To study algebraic structures defined by operations with multiple
inputs and multiple outputs, like bialgebras or Lie bialgebras for
instance, one needs to generalize the notion of operad and
introduce the notion of PROP.

It is natural to try to generalize the Koszul duality for PROPs. A
first result in the direction is due to W. L. Gan in \cite{G}, see
also M. Markl and A. A. Voronov in \cite{MV}.\\

\noindent We work over a field $k$ of characteristic $0$. The symmetric group
on $n$ elements is denoted by $\Sy_n$.

\section{\bf PROPs and $\Po$-gebras}

Over a vector space, various algebraic structures can be
considered like algebras, coalgebras, bialgebras. Following J.-P.
Serre in \cite{Serre}, we call \emph{gebra} any one of these
structures. The set $\Po(m,\, n)$ of the operations of $n$ inputs and $m$ outputs
acting
on a gebra $A$ is a module over $\Sy_m$ on the left and over $\Sy_n$ on
the right.
$$ \Po(m,\,n)\otimes_{\Sy_n} A^{\otimes n} \longrightarrow A^{\otimes m}.$$

\begin{dei}[$\Sy$-bimodule]
An \emph{$\Sy$-bimodule} $(\Po(m,\, n))_{m,\, n \in \mathbb{N}^*}$
is a collection of $\Sy_m\times \Sy_n^{op}$-modules.
\end{dei}

\subsection{The composition product $\boxtimes$}

We introduce a product on $\Sy$-bimodules which describes the composition of
operations.\\

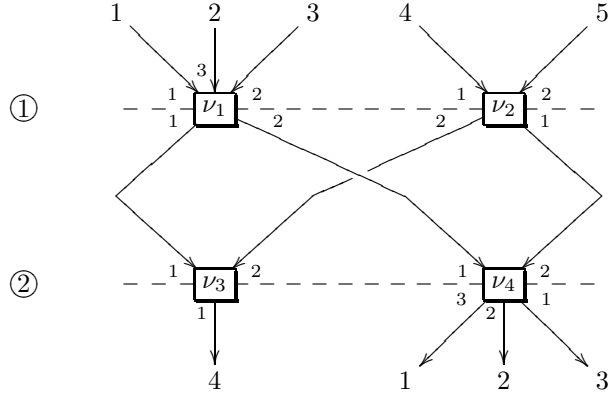
\begin{figure}[h]
$$ \xymatrix{
 & \ar[dr]_(0.7){1} 1&2 \ar[d]_(0.6){3} &3\ar[dl]^(0.7){2} &4\ar[dr]_(0.7){1} & &5\ar[dl]^(0.7){2} \\
*+[o][F-]{1} & \ar@{--}[r]& *+[F-,]{\nu_1} \ar@{-}[dl]_(0.3){1}
\ar@{-}[drr]^(0.3){2} \ar@{--}[rrr]& & & *+[F-,]{\nu_2}
\ar@{-}[dr]^(0.3){1}   \ar@{-}[dll]_(0.3){2} |(0.75) \hole  \ar@{--}[r]& \\
 & *=0{} \ar[dr]_(0.7){1} & &*=0{}\ar[dl]^(0.7){2} &*=0{}\ar[dr]_(0.7){1} & &*=0{} \ar[dl]^(0.7){2}\\
*+[o][F-]{2}& \ar@{--}[r] & *+[F-,]{\nu_3} \ar[d]_(0.3){1}
\ar@{--}[rrr]& & & *+[F-,]{\nu_4} \ar[dl]_(0.3){3} \ar[d]_(0.3){2}
\ar[dr]^(0.3){1}
\ar@{--}[r] & \\
& &4 & &1 &2 &3 } $$ \caption{Example of a 2-levels graph.}
\label{figrapheniveau}
\end{figure}

Let $\G$ be the set of finite graphs with a flow. We suppose that
the inputs and the outputs of each vertex are labeled by integers.
When the vertices of a graph $g$ can be dispatched on two levels,
we denote $\N_i$ ($i=1,\, 2$) the set of vertices belonging to the
$i^{\rm th}$ level. We denote by $\G^2$ the set of such graphs
(\emph{cf.} Figure~\ref{figrapheniveau}).

\begin{dei}[Product $\boxtimes$]
Given two $\Sy$-bimodules $\Po$ and $\Q$, we define their product
by the formula
$$\Po \boxtimes \Q = \left( \bigoplus_{g\in \mathcal{G}^2}
\bigotimes_{\nu \in \N_2} \Po(|Out(\nu)|,\, |In(\nu)|) \otimes
\bigotimes_{\nu \in \N_1}
\Q(|Out(\nu)|,\, |In(\nu)|)
 \right)\Bigg/ \approx,$$
 where the relation $\approx$ is generated by
$$
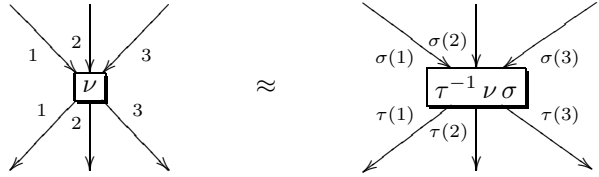
\xymatrix{\ar[dr]_(0.5){1} &\ar[d]_(0.5){2} & \ar[dl]^(0.5){3}
& & \ar[dr]_(0.5){\sigma(1)}& \ar[d]_(0.5){\sigma(2)}&
\ar[dl]^(0.5)
{\sigma(3)} \\
 & *+[F-,]{\nu} \ar[dl]_(0.4){1}\ar[d]_(0.4){2}\ar[dr]^(0.4){3}&
 &\approx & &
 *+[F-,]{\tau^{-1}\, \nu\,\sigma} \ar[dl]_(0.5){\tau(1)}\ar[d]_(0.5){\tau(2)}
 \ar[dr]^(0.5){\tau(3)}&\\
 & & & & & &\ .}$$
Here $|Out(\nu)|$ and $|In(\nu)|$ are the numbers of the outgoing
and the incoming edges of the vertex $\nu$.
\end{dei}

This product has an algebraic writing using the symmetric groups
(\emph{cf.} \cite{V}).

\subsection{PROPs}

The notion of PROP models the operations acting on a certain type
of gebras and their compositions.\\

We denote by $\widetilde{I}$ the identity $\Sy$-bimodule defined
by the formula
$$\left\{
\begin{array}{l}
\widetilde{I}(n,\, n)=k[\Sy_n], \\
\widetilde{I}(m,\, n)=0 \quad \textrm{elsewhere}.
\end{array}
\right.$$

\begin{dei}[PROP]
A structure of \emph{PROP} over an $\Sy$-bimodule $\Po$ is given
by the following data

\begin{itemize}

\item an associative \emph{composition} $\Po\boxtimes \Po
\xrightarrow{\mu} \Po$,

\item a \emph{unit} $\widetilde{I}\xrightarrow{\eta}\Po$.
\end{itemize}
\end{dei}

\begin{Rq}
This definition of a PROP is equivalent to the definition given by
Mac Lane (\emph{cf.} \cite{MacLane})
\end{Rq}

\begin{exem} For any vector space $V$, the sets $(Hom(V^{\otimes n},\,
V^{\otimes m}) )_{m,\,n \in \mathbb{N}^*}$ of morphisms from
$V^{\otimes n}$ to $V^{\otimes m}$ with the composition of
morphisms as $\mu$ form a PROP, denoted $End(V)$.

 The associative algebras
and the operads are examples of PROPs.
\end{exem}

Dually, we define the notion of \emph{coPROP}, which is a PROP in
the opposite category.

\subsection{Quadratic PROPs}

We give a categorical construction of the free monoid in \cite{V},
which applied here, gives the free PROP.

\begin{pro}
The free PROP over an $\Sy$-bimodule $V$, denoted $\F(V)$, is
given by the direct sum on the set of graphs (without
level), where each vertex is indexed by an element of $V$ :

$$\F(V)=\left( \bigoplus_{g\in \mathcal{G}} \bigotimes_{\nu \in
\mathcal{N}} V(|Out(\nu)|,\, |In(\nu)|) \right)\Bigg/\approx.$$
\end{pro}

\begin{Rq} This construction is analytic in $V$. The part of
weight $n$ of $\F(V)$, denoted $\F_{(n)}(V)$, is the direct
summand generated by the graphs with $n$ vertices.
\end{Rq}

Dually, we define the \emph{cofree connected coPROP $\F^c(V)$ on
$V$} with the same underlying $\Sy$-bimodule $\F(V)$.

\begin{dei}[Quadratic PROP]
A \emph{quadratic PROP} $\Po$ is a PROP $\Po=\F(V)/(R)$
generated by an $\Sy$-bimodule $V$ and a space of relations
$R\subset {\F_c}_{(2)}(V)$, where ${\F_c}_{(2)}(V)$ is the direct summand of
$\F(V)$ generated by the connected graphs with $2$ vertices.
\end{dei}

Since the relations $R$ of a quadratic PROP are homogenous,
a quadratic PROP is weight-graded.

\begin{exem}
The PROP of Lie bialgebras, denoted $\BLi$ (\emph{cf.} V. Drinfeld
\cite{D}), the PROP of combinatorial Lie bialgebras, denoted
$\BLi_0$ (\emph{cf.} M. Chas \cite{C})  and the PROP of
infinitesimal Hopf bialgebras, denoted $\IBi$ (\emph{cf.} M.
Aguiar \cite{A}), are examples of quadratic PROPs.
\end {exem}

\subsection{$\Po$-gebras}
The notion of a gebra over a PROP $\Po$ is the generalisation of the notion of
an algebra over an operad.

\begin{dei}[$\Po$-gebra]
A structure of \emph{$\Po$-gebra} over a vector space $A$ is given by a morphism
of PROPs : $\Po \to End(A)$.
\end{dei}

\begin{exem}
A $\BLi$-gebra is exactly a Lie bialgebra, a $\BLi_0$-gebra is a
combinatorial Lie bialgebra and an $\IBi$-gebra is an
infinitesimal Hopf bialgebra.
\end{exem}

\section{\bf Koszul duality}

We generalize the Koszul duality theory of associative algebras and operads to PROPs.

\subsection{Bar and cobar constructions}
We generalize the bar and the cobar constructions of the algebras
and operads to PROPs.

\begin{dei}[Partial product]
For an augmented PROP $(\Po=I\oplus \oPo,\, \mu,\, \eta)$, we
define the \emph{partial product} as the restriction of the
composition $\mu$ to the sub-module of $\Po \boxtimes \Po$
composed by connected graphs with only one vertex on each level
indexed by an element of $\bar{\Po}$ (and the other vertices
indexed by $I$, the unit).
\end{dei}

We denote by $\Sigma$ the homological \emph{suspension}. For
instance, we have $(\Sigma\Po)_{n+1}=\Po_n$, where $n$ is the
homological degree.

\begin{pro}
There exists a unique coderivation $d_\theta$ on $\F^c(\Sigma
\oPo)$ whose restriction on $\F^c_{(2)}(\Sigma \oPo)$ is the
partial product.
\end{pro}

\begin{Rq}
This notion generalizes the \emph{edge contraction}, given by M.
Kontsevich \cite{K}, which defines the boundary of the \emph{graph homology}.
\end{Rq}

\begin{dei}[Bar construction]
Let $(\Po,\, \delta)$ be a differential augmented PROP. The
\emph{bar construction} of $\Po$ is the following chain complex :
$$\B(\Po)=\left( \F^c(\Sigma \oPo),\, \delta + d_\theta \right) .$$
\end{dei}

Dually, we define the \emph{cobar construction} of a differential
co-augmented coPROP $(\Co, \delta)$ and we denote it
$\B^c(\mathcal{C})$.

\subsection{Koszul dual and Koszul complex}
We give the basic definitions of the objects involved in the
Koszul duality theory for PROPs.\\

A PROP in the category of weight-graded vector spaces is called a
\emph{weight-graded PROP}. Quadratic PROPs are examples of
weight-graded PROPS.

\begin{dei}[Koszul dual]
To a weight-graded augmented PROP $\Po$, we associate a dual
coPROP, denoted $\Po^{\ac}$, which is a sub-coPROP of the bar
construction $\Po^{\ac} \hookrightarrow \B(\Po)$.
\end{dei}

\begin{Rq}
Under finite dimensional
hypothesis, the linear dual of $\Po^{\ac}$ gives a PROP which corresponds, in the
cases of associative algebras and operads,
to the classical Koszul dual $\Po^!$ (\emph{cf.} S. Priddy \cite{P} and
V. Ginzburg and M.M. Kapranov \cite{GK}).
\end{Rq}

\begin{dei}[Koszul complex]
On the $\Sy$-bimodule $\Po^{\ac}\boxtimes \Po$, there is defined a
boundary map $d_K$. This chain complex is called the \emph{Koszul
complex} of $\Po$.
\end{dei}

\subsection{Koszul criterion}

We prove a criterion that
determines whether the cobar construction on the dual coPROP
gives a resolution of $\Po$ or not.

\begin{thm}[Koszul criterion]
Let $\Po$ be a differential weight-graded augmented PROP (for
instance a quadratic PROP), the following assertions are
equivalent

\begin{tabular}{ll}

$(1)$ & The Koszul complex $(\Po^{\ac} \boxtimes\Po,\, d_K)$ is acyclic.\\

$(1')$ & The Koszul complex $(\Po\boxtimes \Po^{\ac},\, d'_K)$ is acyclic.\\

$(2)$ & The inclusion $\Po^{\ac} \hookrightarrow \B(\Po)$
is a quasi-isomorphism of differential PROPs.\\

$(2')$ & The projection $\B^c(\Po^{\ac}) \twoheadrightarrow \Po$
is a quasi-isomorphism of differential PROPs.
\end{tabular}
\end{thm}

\begin{sketch}
\hfill

\begin{itemize}

\item We remark that the product $\boxtimes$ induces functors
$A\boxtimes -$ and $-\boxtimes A$ which are analytic functors.
With this graduation and the weight graduation of the PROP $\Po$,
we generalize homological lemmas, called \emph{comparison lemmas},
proved by B. Fresse in \cite{Fresse} for operads, to PROPs.

\item We prove that the augmented bar construction
$\B(\Po)\boxtimes \Po$ of a PROP $\Po$ and the co-augmented cobar
construction $\mathcal{C}\boxtimes B^c(\mathcal{C})$ are acyclic.

\item We apply the comparison lemmas to show that the bar-cobar
construction over a weight-graded PROP $\Po$ gives a resolution of
$\Po$.
$$\B^c(\B(\Po))\longrightarrow \Po.$$

\item We simplify the bar-cobar construction with the comparison lemmas to conclude.

\end{itemize}
 $\fin$
\end{sketch}

\begin{Rq}
This theorem includes the cases of associative algebras and operads.\\
\end{Rq}

A PROP $\Po$ that verifies these assertions is called a
\emph{Koszul PROP}. In this case, the cobar construction on the
dual $\Po^{\ac}$ is the resolution of $\Po$ called the
\emph{minimal model} of $\Po$.

\begin{exem}
We prove that the PROPs $\BLi$, $\BLi_0$ and $\IBi$ are Koszul
PROPs. This result can be interpreted in terms of graph homology
like in \cite{MV}.
\end{exem}

\subsection{$\Po$-gebras up to homotopy}

One of the main application of the Koszul duality for PROPs is the definition of a
$\Po$-gebra up to homotopy.

\begin{dei}[$\Po$-gebra up to homotopy]
Let $\Po$ be a Koszul PROP. A gebra over the cobar construction
$\B(\Po)$ of $\Po$ is called a \emph{$\Po$-gebra up to homotopy}
and denoted a \emph{$\Po_\infty$- gebra}.
\end{dei}

\begin{exem}
Applied to the examples given below this defines the notions of
Lie bialgebras, combinatorial Lie bialgebras and infinitesimal
Hopf bialgebras up to homotopy.
\end{exem}


\end{document}